\documentclass[12pt]{article}
\usepackage{hyperref}
\usepackage{amsfonts}
\usepackage{enumerate}
\usepackage{graphicx}
\usepackage{amsmath}
\begin{document}

\begin{center}
{\large\bf On stability of generalised systems of difference equation with non-consistent initial conditions}

\vskip.20in
Nicholas Apostolopoulos$^{1}$, Fernando Ortega$^{2}$ and\ Grigoris Kalogeropoulos$^{3}$\\[2mm]
{\footnotesize
$^{1}$National Technical University of Athens, Greece\\
$^{2}$ Universitat Autonoma de Barcelona, Spain\\
$^{3}$National and Kapodistrian University of Athens, Greece}
\end{center}

{\footnotesize
\noindent
\textbf{Abstract:}
For given non-consistent initial conditions, we study the stability of a class of generalised linear systems of difference equations with constant coefficients and taking into account that the leading coefficient can be a singular matrix. We focus on the optimal solutions of the system and derive easily testable conditions for stability. 
\\\\[3pt]
\\[3pt]

\vskip.2in

\section{Introduction}
Singular systems of difference/differential equation have been studied by many authors in the past years. See [1-28], and [29-36] for recent applications of such systems. For an extended version of this type of systems using fractional operators, see [37-45]. In this article, we consider the following initial value problem:
\begin{equation}\label{eq1}
\begin{array}{cc}
FY_{k+1}=GY_k, & k= 1, 2,...,\\\\
Y_{0}.
\end{array}
\end{equation}      
Where  $F, G \in \mathbb{R}^{m \times m}$ and $Y_k\in \mathbb{R}^{m}$. The matrix $F$ is singular (det$F$=0). The initial conditions $Y_0$ are considered to be non-consistent. Note that the initial conditions are called consistent if there exists a solution for the system which satisfies the given conditions. We will also assume that the pencil of the system is regular, i.e. for an arbitrary $s\in\mathbb{C}$ we have det$(sF-G)\neq 0$, see [46-53].

There are stability results in the literature dealing with regular systems and for generalised systems with consistent conditions, see [11, 12, 47, 48]. As already mentioned, we consider initial conditions that are non-consistent. This means that the system has infinite many solutions and an optimal solution is required for this case, see [54, 55]. The aim of this paper is to study the stability of the optimal solution of \eqref{eq1}.

\section{Preliminaries}
Some tools from matrix pencil theory will be used throughout the paper. Since in this article we consider the system \eqref{eq1} with a \textsl{regular pencil}, the class of $sF-G$ is characterized by a uniquely defined element, known as the Weierstrass canonical form, see [46-53], specified by the complete set of invariants of $sF-G$. This is the set of elementary divisors of type  $(s-a_j)^{p_j}$, called \emph{finite elementary divisors}, where $a_j$ is a finite eigenvalue of algebraic multiplicity $p_j$ ($1\leq j \leq \nu$), and the set of elementary divisors of type $\hat{s}^q=\frac{1}{s^q}$, called \emph{infinite elementary divisors}, where $q$ is the algebraic multiplicity of the infinite eigenvalue. $\sum_{j =1}^\nu p_j  = p$ and $p+q=m$.
\\\\
From the regularity of $sF-G$, there exist non-singular matrices $P$, $Q$ $\in \mathbb{R}^{m \times m}$ such that 
\begin{equation}\label{eq3}
\begin{array}{c}PFQ=\left[\begin{array}{cc} I_p&0_{p,q}\\0_{q,p}&H_q\end{array}\right],
\\\\
PGQ=\left[\begin{array}{cc} J_p&0_{p,q}\\0_{q,p}&I_q\end{array}\right].\end{array}
\end{equation}
$J_p$, $H_q$ are appropriate matrices with $H_q$ a nilpotent matrix with index $q_*$, $J_p$ a Jordan matrix and $p+q=m$. With $0_{q,p}$ we denote the zero matrix of $q\times p$. The matrix $Q$ can be written as
\begin{equation}\label{eq4}
Q=\left[\begin{array}{cc}Q_p & Q_q\end{array}\right].
\end{equation}
$Q_p\in \mathbb{R}^{m \times p}$ and $Q_q\in \mathbb{R}^{m \times q}$. 
The following results have been proved.
\\\\
\textbf{Theorem 2.1.}  (See [1-28]) We consider the systems \eqref{eq1} with a regular pencil. Then, its solution exists and for $k\geq 0$, is given by the formula
\begin{equation}\label{eq6}
    Y_k=Q_pJ_p^kC.  
\end{equation}
Where $C\in\mathbb{R}^p$ is a constant vector. The matrices $Q_p$, $Q_q$, $J_p$, $H_q$ are defined by \eqref{eq3}, \eqref{eq4}. 
\\\\

\subsection*{Non-consistent initial conditions}

The following proposition identifies if the initial conditions are non-consistent:
\\\\
\textbf{Proposition 2.1.} The initial conditions of system \eqref{eq1} are non-consistent if and only if 
\[
Y_0\notin colspan Q_p.
\]
We can state the following Theorem, see [54, 55]. 
\\\\
\textbf{Theorem 2.2.} We consider the system \eqref{eq1} with known non-consistent initial conditions. For the case that the pencil $sF-G$ is regular, after a perturbation to the non-consistent initial conditions accordingly
\[
min\left\|Y_0-\hat Y_0\right\|_2,
\]
or, equivalently,
\[
\left\|Y_0-Q_p(Q_p^*Q_p)^{-1}Y_0\right\|_2,
\]
an optimal solution of the initial value problem \eqref{eq1} is given by
\begin{equation}\label{eq7}
    \hat Y_k=Q_pJ_p^k(Q_p^*Q_p)^{-1}Q_p^*Y_0.
\end{equation}
The matrices $Q_p$, $J_p$ are given by by \eqref{eq3}, \eqref{eq4}.

\section{Main Results}

We will focus on the stability of equilibrium state(s) of homogeneous singular discrete time systems:
\\\\
\textbf{Definition 3.1.} For any system of the form \eqref{eq1} , $Y_*$ is an equilibrium state if it does not change under the initial condition, i.e.: $Y_*$ is an equilibrium state if and only if $Y_{0}=Y_*$ implies that $Y_k=Y_*$ for all $k \geq$ 1.
\\\\
The set of equilibrium states for a given singular linear system in the form of \eqref{eq1} are given by the following Proposition, see [11, 12]:
\\\\
\textbf{Proposition 3.1.} Consider the system \eqref{eq1}. Then if 1 is not
an eigenvalue of the pencil $sF-G$ then
\[
Y_*=0_{m,1}
\]
is the unique equilibrium state of the system \eqref{eq1}. If 1 is a finite
eigenvalue of $sF-G$, then the set $E$ of the
equilibrium points of the system \eqref{eq1} is the vector space defined by
\[
E=N_r(F-G)\cap {\rm colspan} Q_p.
\]
Where $N_r$ is the right null space of the matrix $F-G$, $Q_p$ is a
matrix with columns the $p$ linear independent (generalized)
eigenvectors of the $p$ finite eigenvalues of the pencil.
\\\\
\textbf{Proof.} If $Y_*$ is an equilibrium state of system \eqref{eq1}, then this implies that for
\[
Y_0=Y_*
\]
we have
\[
Y_*=Y_k=Y_{k+1}.
\]
If 1 is not an eigenvalue of the pencil then det$(F-G)\neq0$ and 
\[
FY_*=GY_*,
\]
or, equivalently,
\[
(F-G)Y_*=0_{m,1}.
\]
Then the above algebraic system has the unique solution
\[
Y_*=0_{m,1}.
\]
which is the unique equilibrium state of the system. If 1 is a
finite eigenvalue of the pencil then det$(F-G)$=0. If
$Y_*$ is an equilibrium state of the system, then this implies that
for
\[
Y_0=Y_*
\]
we have
\[
Y_*=Y_k=Y_{k+1}.
\]
This requires that $Y_*$ must be a consistent initial condition which from Proposition 2.1 is equal to
\[
Y_*\in colspan Q_p.
\]
Moreover we have
\[
FY_*=GY_*,
\]
or, equivalently,
\[
(F-G)Y_*=0_{m,1},
\]
or, equivalently,
\[
Y_*\in N_r(F-G).
\]
Hence,
\[
Y_*\subseteq N_r(F-G)\cap colspan Q_p,
\]
or, equivalently,
\[
E\subseteq N_r(F-G)\cap {\rm colspan} Q_p.
\]
Let now $Y_*\in N_r(F-G)\cap colspan Q_p$ then we can consider
\[
Y_0=Y_*
\]
as a consistent initial condition and
\[
(F-G)Y_*=0_{m,1}
\]
or, equivalently,
\[
FY_*=GY_*,
\]
where $Y_*$ is solution of the system and combined with $Y_0=Y_*$ we have $Y_*\in E$, or, equivalently,
\[
N_r(F-G)\cap {\rm colspan} Q_p\subseteq E
\]
The proof is completed.
\\\\
\textbf{Theorem 3.1.} We consider the system \eqref{eq1} with non-consistent initial conditions. An optimal solution is then given by \eqref{eq7}.
Then an equilibrium state $Y_*\in  E$ is stable in the sense of Lyapounov, if and
only if, there exist a constant $c \in (0, +\infty)$, such that
$\left\|J_p^k\right\| \leq c < +\infty$, for all $k\geq 0$.
\\\\
\textbf{Proof.} An optimal solution is then given by \eqref{eq7}:
\[
\hat Y_k =
 Q_pJ_p^k(Q_p^*Q_p)^{-1}Q_p^*Y_0.
 \]
 We assume that there exist a constant $c \in (0, +\infty)$ such that $\left\|J_p^k\right\|$ $\leq c < +\infty$, for all $k> 0$. Furthermore let an equilibrium state $Y_* \in E$. Then
 \[
 Y_*=Q_pJ_p^k(Q_p^*Q_p)^{-1}Q_p^*Y_*
 \]
 and easy we obtain
 \[
 \hat Y_k-Y_*=Q_pJ_p^k(Q_p^*Q_p)^{-1}Q_p^*(\hat Y_k-Y_*),
 \]
or, equivalently,
\[
\hat Y_k-Y_*=\left[\begin{array}{cc}Q_p&0_{m,q}\end{array}\right]\left[\begin{array}{cc}J_p^k&0_{p,q}\\0_{q,p}&0_{q,q}\end{array}\right]\left[\begin{array}{c}(Q_p^*Q_p)^{-1}Q_p^*\\0_{q,m}\end{array}\right](\hat Y_k-Y_*)
\]
If we set $\left\|Q_p\right\|=\left\|\left[\begin{array}{cc}Q_p&0_{m,q}\end{array}\right]\right\|$ and $\left\|(Q_p^*Q_p)^{-1}Q_p^*\right\|=\left\|\left[\begin{array}{c}(Q_p^*Q_p)^{-1}Q_p^*\\0_{q,m}\end{array}\right]\right\|$. Then by taking norms for every k$\geq 0$ we have
 \[
 \left\|\hat Y_k-Y_*\right\|\leq\left\|Q_p\right\|\left\|J_p^k\right\|\left\|(Q_p^*Q_p)^{-1}Q_p^*\right\|\left\|Y_0-Y_*\right\|
 \]
 Hence for any $\epsilon>0$, if we chose $\delta(\epsilon)$=$\frac{\epsilon}{\left\|Q_p\right\|\left\|(Q_p^*Q_p)^{-1}Q_p^*\right\|c}$, then 
 \[
 \left\|Y_0-Y_*\right\|\leq \delta(\epsilon)
 \]
 implies that for every $\epsilon>0$
 \[
 \left\|\hat Y_k-Y_*\right\|\leq\left\|Q_p\right\|\left\|J_p^k\right\|\left\|(Q_p^*Q_p)^{-1}Q_p^*\right\|\left\|Y_0-Y_*\right\|,
 \]
 or, equivalently,
  \[
 \left\|\hat Y_k-Y_*\right\|\leq\left\|Q_p\right\|c\left\|(Q_p^*Q_p)^{-1}Q_p^*\right\|\frac{\epsilon}{\left\|Q_p\right\|\left\|(Q_p^*Q_p)^{-1}Q_p^*\right\|c}\leq \epsilon,
 \]
 or, equivalently,
 \[
 \left\|\hat Y_k-Y_*\right\|\leq \epsilon.
 \]
The proof is completed.
\\\\
\textbf{Theorem 3.2.} The system \eqref{eq1} with non-consistent initial conditions has an optimal solution is then given by \eqref{eq7}. Then it is asymptotic stable at large, if and only if, all the finite eigenvalues of $sF-G$ lie within the open disc,
 \[
 \left|s\right|<1.
 \]
 \textbf{Proof.} The system \eqref{eq1} with non-consistent initial conditions has an optimal solution is then given by \eqref{eq7}:
 \[
 \hat Y_k=Q_pJ_p^k(Q_p^*Q_p)^{-1}Q_p^*Y_0.
 \]
 Let $a_j$ be a finite eigenavalue of the pencil with algebraic multiplicity $p_j$. Then the Jordan matrix $J_p^k$ can be written as
 \[
 J_p^k={\rm blockdiag}\left[\begin{array}{cccc} J_{p_1}^{k}(a_1)& J_{p_2}^{k}(a_2)& \dots& J_{p_\nu}^{k}(a_\nu)\end{array}\right],
 \]
 with $J_p^{k}\in\mathcal{M}_{p_j}$ be a Jordan block. Every element of this matrix has the specific form
 \[
 k^{p_j}a_j^{k}.
 \]
 The sequence
 \[
 k^{p_j}\left|a_j^{k}\right|,
 \]
 can be written as
 \[
 k^{p_j}e^{kln\left|a_j\right|}.
 \]
The system is asymptotic stable at large, when
 \[
 \lim _{k\longrightarrow \infty} \hat Y_k=Y_*.
 \]
 Thus this holds if and only if
 \[
 \ln\left|a_j\right| < 0,
 \]
 or, equivalently,
 \[
 \left|a_j\right| < 1.
 \]
 Then for $k\rightarrow+\infty$:
 \[
 k^{p_j}e^{k\left|\ln a_j\right|}\rightarrow 0,
 \]
 or, equivalently,
 \[
 k^{p_j}\left| a_j\right|^{k}\rightarrow 0,
 \]
 or, equivalently, for every $k\geq 0$
 \[
 J_p^{k}\longrightarrow 0_{p, p}.
 \]
 Then for every initial condition $Y_{0}$
 \[
 \lim _{k\longrightarrow \infty} \hat Y_k=0_{m,1}.
 \]
The proof is completed.

\section*{Conclusions}
In this article we focused and provided properties for the stability of the optimal solutions of a linear generalized discrete time system in the form of \eqref{eq1} for given non-consistent initial conditions.

\end{document}